\documentclass[10pt]{amsart}
\usepackage{amssymb, amsmath, amsfonts}
\usepackage{mathtools}
\usepackage{amscd}
\usepackage{verbatim}
\usepackage{color}
\usepackage[color,all]{xy}
\usepackage[all]{xy}
\usepackage{tikz-cd}
\usepackage{kotex}
\numberwithin{equation}{section}

\newcommand{\PP}{\mathbb P}
\newcommand{\pp}{\mathbb P}
\newcommand{\cc}{\mathbb C}

\newcommand{\td}{\tilde{d}}
\newcommand{\ta}{\widetilde{\alpha}}
\newcommand{\tb}{\widetilde{\beta}}

\newcommand{\C}{\mathcal{C}}

\newcommand{\V}{\mathcal V} 

\newcommand{\cL}{\mathcal L}

\newcommand{\Sym}{\mathrm{Sym}}
\newcommand{\rank}{\mathrm{rk}}

\newcommand{\Hom}{\mathrm{Hom}}

 \newcommand{\cM}{\mathcal{M}}
\newcommand{\cV}{\mathcal{V}}

\def\Hom{\mathrm{Hom} }

\newcommand{\M}{\mathcal{M}}

\newcommand{\tE}{\widetilde{E}}

\newcommand{\tV}{\widetilde{V}}
\newcommand{\tom}{\widetilde{\omega}}

\newcommand{\rk}{\mathrm{rk}}

\newcommand{\cO}{\mathcal{O}}

\newcommand{\MS}{\mathcal MS}
\newcommand{\MO}{\mathcal MO}

\newcommand{\cE}{{\mathcal E}}
\newcommand{\cF}{{\mathcal F}}

\usepackage{hyperref}
\usepackage{color}
\usepackage[color,all]{xy}
\hypersetup{colorlinks, citecolor=red, linkcolor=blue}

\newtheorem{theorem}{{\textbf Theorem}}[section]
\newtheorem{prop}[theorem]{{\textbf Proposition}}
\newtheorem{cor}[theorem]{{\textbf Corollary}}
\newtheorem{lemma}[theorem]{{\textbf Lemma}}

\newtheorem{definition}[theorem]{{\textbf Definition}}
\newenvironment{defn}{\begin{definition}\rm}{\end{definition}}

\newtheorem{rmk}[theorem]{{\textbf Remark}}
\newenvironment{remark}{\begin{rmk}\rm}{\end{rmk}}
\newtheorem{assume}[theorem]{{\textbf Assumption}}

\title[Minimal rational cuvres]{Minimal rational curves \\on the moduli spaces of \\symplectic   and orthogonal bundles}
\author{Insong Choe}
\address{Department of Mathematics\\
  Konkuk  University\\
120 Neungdong-ro, Gwangjin-gu, Seoul 05029, Korea}
\email[I.~Choe]{ischoe@konkuk.ac.kr}
\author{Kiryong Chung}
\address{Department of Mathematics Education\\
  Kyungpook National University, 80 Daehakro, Bukgu, Daegu 41566, Korea}
\email[K.~Chung]{krchung@knu.ac.kr}
\author{Sanghyeon Lee}
\address{Department of Mathematics\\
  Korea Institute for Advanced Study (KIAS)\\
  85 Hoegiro, Dongdaemun-gu, Seoul 02455, Korea}
\email[S. ~Lee]{sanghyeon@kias.re.kr}

\begin{document}

\maketitle

\begin{abstract}
Let $C$ be an algebraic curve of genus $g$ and $L$ a line bundle over $C$.
Let $\MS_C(n,L)$ and $\MO_C(n,L)$ be the moduli spaces of $L$-valued symplectic and orthogonal bundles respectively, over $C$ of rank $n$. We construct rational curves of Hecke type on these moduli spaces which generalize the Hecke curves on the moduli space of vector bundles. As a main result, we show that these  curves have the minimal degree among the rational curves passing through a general point of the moduli spaces. As its byproducts, we show the non-abelian Torelli theorem and compute the automorphism group of  the moduli spaces. 
\end{abstract}

\section{Introduction}

Let $C$ be a smooth algebraic curve of genus $g \ge 2$ over $\cc$. For integers $n \ge 2$ and $d$, the moduli space $SU_C(n, d)$ of semistable bundles over $C$ of rank $n $ with  a fixed determinant of degree $d$ is known to be a Fano variety of Picard number one.  Xiaotao Sun \cite{Sun} proved that under the assumption that $g \ge 3$ (except  the case when $g=3, n=2$, and $d$ even), a rational curve passing through a general point of $SU_C(n, d)$ has the minimal degree if and only if it is a so called \emph{Hecke curve}. This  readily yields a simple proof of nonabelian Torelli theorem and the description of the automorphism group  of $SU_C(n, d)$ (\cite[Corollary 1.3 and 1.4]{Sun}).
The Hecke curves have been widely studied for $n=2$, for which case the degree minimality was proven  in \cite[Proposition 8]{Hw1}.  For introduction to Hecke curves for $n >2$, see \cite{Hw} and \cite{HR}. 

The goal of this paper is to establish a similar result for symplectic bundles and  orthogonal bundles. Let us briefly explain the main results, pending the precise definition and  backgrounds for symplectic   and orthogonal bundles  to \S\:\ref{subsecor}.

Let $\MS_C(n, L)$ (resp. $\MO_C(n,L)$) be the moduli space of $L$-valued symplectic  (resp. orthogonal) bundles over $C$ of rank $n$ for  a fixed line bundle $L$ over $C$ of degree $\ell$. We first construct variants of Hecke curves on these moduli spaces, called \emph{symplectic}  and \emph{orthogonal Hecke curves}, which  provide   covering families of rational curves for $\MS_C(n, L)$ and $\MO_C(n,L)$ respectively. Interestingly,  it turns out that the  orthogonal Hecke curves correspond to rational curves having  degree higher than that of the Hecke curves on  $SU_C(n, d)$, while the symplectic Hecke curves   are just special type of Hecke curves on  $SU_C(n, d)$.

 Since the Picard groups of $\MS_C(n,L)$ and $\MO_C(n,L)$ are infinitely cyclic (under the  assumption on $n$ below), the notion of minimality of the degree of a curve  makes sense without a reference  to a specific ample line bundle.  
The main result of this paper  is as follows:

\begin{theorem}  \label{main1}
 Assume $g \ge 3, \: n \ge 4$ in the symplectic  case and  $g \ge 5, \: n \ge 5$  in the orthogonal case.
Among the   rational curves passing through  a general point of $\MS_C(n,L)$ {\rm(}resp. $\MO_C(n,L)${\rm)}, a  curve    has  minimal degree if and only if it is a symplectic {\rm(}resp. orthogonal{\rm)} Hecke curve.
\end{theorem}

More detailed version  of this statement  is given in Theorem \ref{mainHecke} and Theorem \ref{mainHeckeortho} for symplectic  and orthogonal cases, respectively.

The reason here for the assumption on  the rank $n$   is as follows: 
A symplectic bundle must have  even rank, and $\MS_C(2,L) \cong SU_C(2, L)$  since every vector bundle $V$ of rank 2 has a  symplectic structure coming from the isomorphism $V \cong V^* \otimes (\det V)$.  
Also  orthogonal bundles of rank $\le 4$ are constructed from line bundles and vector bundles of rank 2 (see \cite[p.185]{Mum} and \cite[\S\:4, 5]{Ser2}).   

The idea of proof of Theorem \ref{main1} is as follows.
Once the relevant variants of Hecke curves are defined, we can follow   \cite[\S\:2]{Sun} to compute  the degree of rational curves.
One additional  ingredient is the  Harder--Narasimhan filtration of symplectic and orthogonal bundles. In particular, we observe that   a symplectic or orthogonal bundle over $\pp^1$ must have a special kind of  splitting type. This enables us to adapt the formulas in \cite{Sun}, originally for vector bundles, to the context of symplectic and orthogonal bundles.

To guarantee that the involved maps from $\pp^1$ to the moduli spaces are (generically) injective,  we need to handle a  technical issue regarding  generic isomorphisms between symplectic and orthogonal bundles.  This can be done by a rather long but elementary argument based on dimension counting. 

Theorem \ref{main1} suggests that the symplectic and orthogonal Hecke curves play the same role as the Hecke curves on $SU_C(n,d)$ as in \cite{Hw}, \cite{HR} and \cite{Sun}. 
As  samples,  we will discuss two  consequences  in \S\:\ref{applications}:  NonabelianTorelli theorem (Theorem \ref{Torelli})  and the description of automorphisms of the moduli space (Theorem \ref{automorphism}).   These results are not new, but our approach provides another geometric view to solve these problems.

This paper is organized as follows. In \S\:\ref{subsecor}, the basic facts on symplectic and orthogonal bundles will be explained and their moduli spaces are described.  In  \S\:\ref{Harder}, the description of the symplectic   and orthogonal version of Harder-Narasimhan filtrations  will be given in a concrete way.
In \S\:\ref{Heckesection},  we  construct relevant variants of Hecke curves. We prove Theorem \ref{main1} in \S\:\ref{mainsection} and discuss its applications in \S\:\ref{applications}.

\section{Moduli spaces of symplectic and orthogonal  bundles}\label{subsecor}

Let $L$ be a line bundle over a smooth algebraic curve $C$ of genus $g \ge 0$.
\begin{defn}
A vector bundle $V$ over $C$ of rank $n$  is  called an $L$-valued \textit{symplectic bundle} (resp. \textit{orthogonal bundle}) if it is equipped with a bundle map  $\omega \colon V\otimes V \to L$ whose restriction to each fiber is a
non-degenerate skew-symmetric (resp. symmetric) bilinear form.
\end{defn}
 An $L$-valued symplectic bundle (resp. orthogonal bundle) can also be viewed as a principal $\mathrm{Gp}(n, \cc)$-bundle (resp. $\mathrm{GO}(n, \cc)$-bundle); see \cite{BG1} and \cite{BG2}.
From the isomorphism $V \cong V^* \otimes L$, we have $\det (V)^2 \cong L^{n}$ and so $\deg (V) = \frac{1}{2} n \ell$, where $\ell = \deg (L)$. 

\begin{defn} A subbundle $E$ of $V$ is said to be \textit{isotropic} if $\omega|_{E \otimes E} \equiv 0$.\end{defn} By  linear algebra, the rank of an isotropic subbundle of $V$ cannot exceed the half of $\mathrm{rk}(V)$.  
\begin{defn}
As an symplectic or orthogonal   bundle,  $V$ is \textit{stable} (resp. \textit{semistable}) if $\mu (V) > \mu (E)$ (resp. $\mu (V) \ge \mu (E)$) for every nonzero isotropic subbundle $E$, where $\mu$ is the \emph{slope} of the vector bundle defined by $\mu(V) = \frac{ \deg(V)}{\rk (V)}$. 
 \end{defn}

When $g \ge 2$, the semistable $L$-valued symplectic (resp. orthogonal) bundles over $C$ of rank $n$ form a moduli space denoted by $\MS_C(n, L)$ (resp. $\MO_C(n, L)$).

Note  that every vector bundle $V$ of rank 2 has a $\det(V)$-valued symplectic structure coming from the isomorphism $V \cong V^* \otimes \det (V)$. Hence  we have $\MS_C(2,L) \cong SU_C (2, \ell)$.
In general, the value of the  symplectic form $\omega \colon \wedge^2 V \to L$  determines  $\det (V)$  via the  isomorphism given by the composition: 
\[
\det (V) = \wedge^nV \hookrightarrow (\wedge^2 V)^{\otimes \frac{n}{2} }  \xrightarrow{\omega \otimes \cdots \otimes \omega  } L^{\frac{n}{2} }.
\]
 It is known that 
 $\MS_C(n, L)$ is irreducible of dimension $  \frac{1}{2} n(n+1)(g-1)$.
  
On the other hand, $\MO_C(n, L)$ has several irreducible components:  An $L$-valued orthogonal bundle $V$ has a determinant  such that $\det (V)^2 \cong L^{n}$, hence  there are many  components  distinguished by  $c_1(V)$. Also for a fixed $c_1(V)$,  there is another  topological invariant called the second Stiefel--Whitney class $w_2(V) \in \mathbb{Z}_2$ 
(see \cite{Ser}  for details).  It is known that all the components of 
 $\MS_C(n, L)$ have dimension $  \frac{1}{2} n(n-1)(g-1)$. For dimension of the moduli spaces, see \cite[Theorem 5.9]{Ramanathan}.

Since most of the arguments in this paper will go parallel for the symplectic and orthogonal   case, we simply write $\M$ to denote  either $\MS_C(n, L)$ or an irreducible component of $\MO_C(n, L)$, unless there is   a possibility of confusion.


The notion of the minimality of the degree of a curve on $\M$  is independent of the choice of an ample line bundle on $\M$, due to the following fact: 
\begin{lemma}  Assume $n \ge 4$ for symplectic case and $n \ge 5$ for orthogonal case.
Then the Picard group of $\M$ is infinitely cyclic. 
\end{lemma}
\begin{proof}
This is shown in \cite{BLS}.  We remark that this does not hold in the orthogonal case for $n=4$  (see \cite[\S\:4, 5]{Ser2}). \end{proof}

 In this paper, we choose the reference  line  bundle   as follows. Let
\[
f \colon \M \to SU_C(n, {\textstyle \frac{1}{2}} n \ell) 
\]
be the forgetful morphism sending an $L$-valued symplectic or orthogonal bundle to the underlying vector bundle of degree $ { \frac{1}{2}} n \ell$, where $\ell = \deg (L)$.  
A priori, $f$ is a rational map defined only on the points represented by symplectic or orthogonal bundles  whose underlying vector bundles are semistable.  But   it is known that  every semistable symplectic or orthogonal bundle is semistable as a vector bundle (see Lemma \ref{stablesemistable}).  Furthermore, $f$ is a generically injective morphism, since a  simple vector bundle $W$ cannot have more than one $L$-valued symplectic or orthogonal structure (giving an isomorphism    $W \cong W^* \otimes L$). 
 
 Now the reference line bundle will be the pull-back of the anti-canonical line bundle on $SU_C(n, {\textstyle \frac{1}{2}} n \ell) $. As an advantage of this choice, we can compare the degree of rational curves on the moduli spaces $\M$ and $SU_C(n, {\textstyle \frac{1}{2}} n \ell) $.   
\begin{defn}
Let $u \colon \pp^1 \to \M$ be a generically injective map which gives a rational curve  passing through a general point of $\M$. 
We define its \textit{degree}  as the degree of $(f \circ u)^* (-K)$, where  $K$ is the canonical line bundle  on $SU_C(n,  \frac{1}{2}n \ell) $. 
\end{defn}

\section{Harder--Narasimhan filtrations} \label{Harder}
\noindent In this section, we describe the Harder--Narasimhan filtration of symplectic and  orthogonal bundles. Most of the results in this section  can be found in the literature, written  in the language of principal bundles in general. For our convenience later, we rewrite  the results in the language of symplectic and  orthogonal  bundles.  
\textbf{The discussion in this section applies to any genus $g \ge 0$ and any rank $n \ge 2$.}

First we recall: 
\begin{lemma} \label{stablesemistable}
Let $V$ be an $L$-valued  symplectic/orthogonal bundle with $\deg V = \frac{1}{2} n \ell$, where $\ell = \deg L$. 
\begin{enumerate}
\item A symplectic/orthogonal bundle $V$ is semistable  if and only if it is semistable as a vector bundle.
\item If a symplectic/orthogonal bundle $V$ is stable, then  the underlying vector bundle is polystable: a direct sum of  stable bundles of the same slope.
\end{enumerate}
\end{lemma}
\noindent Hence to show the semistability of $V$, it suffices to check that there is no destabilizing {isotropic} subbundles.\\

\begin{proof} 
The argument below is borrowed from \cite[proposition 4.2]{Ra}, in which the orthogonal case was discussed.
Assume $V$ is semistable as a symplectic/orthogonal bundle. Let $E$ be any subbundle of rank $r$.
From  the isomorphism $V/E \cong (E^\perp)^* \otimes L$, we observe:
\begin{equation} \label{degeq}
\deg (E^\perp) = \deg (E) +(\frac{n}{2}-r) \ell.
\end{equation}
 Let $M$ and $N$  be the subbundles generated by $E \cap E^\perp$ and $E + E^\perp$, respectively.
If $M =0$, then the map $E \oplus E^\perp \to V$ gives a subsheaf  whose quotient is a torsion sheaf. Hence 
\[
2 \deg E + (\frac{n}{2}-r) \ell \: \le \: \deg (V) = \frac{1}{2}n \ell,
\]
and so $\mu (E)  = \frac{\deg (E)}{r} \le \frac{\ell}{2} = \mu (V)$.

  If $M$ has rank $k >0$, 
 there is an exact sequence
\[
0 \to M \to E \oplus E^\perp \to N \to 0.
\]
 Since the symplectic or orthogonal pairing $ \langle E+ E^\perp, E \cap E^\perp \rangle $ is zero, we have $M \subset N^\perp$. Since $M$ and $N^\perp$ have the same rank, they must coincide. Therefore, we have
\[
\rk( E) + \rk (E^\perp) =  \rk (M) +  \rk (N) =  \rk (M) + \rk ( M^\perp ).
\]
Applying (\ref{degeq}) to both $E$ and $M$, we have
\begin{equation} \label{degM} 
\deg (E) = \deg( M) +  \frac{1}{2} (r-k)  \ell.
\end{equation}
Since $M = E \cap E^\perp$ is isotropic, 
\begin{equation} \label{semistableineq}
\deg (M) \le k \mu(V) = \frac{1}{2} k \ell.
\end{equation}
 Hence
$\deg (E) \le \frac{1}{2} r \ell$ and  $V$ is semistable as a vector bundle.  
 
 When we assume $V$ is a stable symplectic/orthogonal bundle, (\ref{semistableineq}) becomes strict inequality. Hence  if the equality $\mu (E) = \mu(V)$ holds, then $M=0$ and in this case the map  
 $E \oplus E^\perp \to V$ is an isomorphism by \eqref{degeq}. This shows the second claim.
\end{proof}
The following gives the Harder--Narasimhan filtration of a symplectic/orthogonal bundle.
\begin{prop} \label{HN}
 For any  $L$-valued symplectic (resp. orthogonal) bundle $V$, there is a chain of isotropic subbundles:
\begin{equation} \label{seq}
0 =E_0 \subset E_1 \subset E_2 \subset \cdots \subset E_k \subset  V
\end{equation}
such that  
\begin{equation} \label{slopes}
\mu(V / E_k) < \mu(E_k/E_{k-1} ) < \cdots < \mu (E_2/E_1) < \mu(E_1).
\end{equation}
Moreover,  $(E_i/E_{i-1}) \oplus (E_{i-1}^\perp / E_i^\perp)$ for $1 \le i \le k$ and  $E_k^\perp / E_k$ are semistable $L$-valued symplectic (resp. orthogonal) bundles under the inherited symplectic (resp. orthogonal) structure.

As a consequence, the Harder--Narasimhan filtration of the underlying vector bundle of $V$ is given by: 
\begin{equation} \label{HNvb}
0 =E_0 \subsetneq E_1 \subsetneq E_2 \subsetneq \cdots \subsetneq E_k \subseteq E_k^\perp \subsetneq E_{k-1}^\perp \subsetneq \cdots \subsetneq E_1^\perp \subsetneq E_{0}^\perp = V,
\end{equation} 
where the middle term $E_k$ may or may not be equal to $E_k^\perp$.
\end{prop}
\begin{proof} 
If $V$ itself is semistable, there is nothing to prove.
If $V$ is  unstable,  let $E_1$ be an  isotropic subbundle of $V$ such that
\begin{itemize}
\item $\mu(E_1)$ is maximal among the isotropic  subbundles  and
\item $\rk (E_1)$ is  maximal  among them.
\end{itemize}
Since $E_1$ has maximal slope and every subbundle of $E_1$ is isotropic, it must be a semistable vector bundle.  

We stop here if $V$ has no destabilizing isotropic subbundle $S$ containing $E_1$. Otherwise, let $E_2$ be an isotropic subbundle containing $E_1$ such that
\begin{itemize}
\item $\mu(E_2)$ is maximal among the isotropic subbundles containing $E_1$ and
\item $\rk (E_2)$ is  maximal  among them.
\end{itemize}
By the maximality of slope, the quotient $E_2 / E_1$ is a semistable vector bundle. Also, since $\mu(E_2) < \mu(E_1)$, we have
$\mu(E_2/E_1) < \mu(E_1).$\footnote{One can use the convenient fact that given a subbundle $A \subset B$, the values of $\mu(B)-\mu(A)$, $\mu(B/A) - \mu(B)$, and $\mu(B/A) - \mu(A)$ have the   same sign.}

Repeating this process, we have the chain
\begin{equation} \label{chain}
0=E_0 \subset E_1 \subset E_2 \subset \cdots \subset E_k \subset V
\end{equation}
of isotropic subbundles  satisfying (\ref{slopes}) such that each $E_i/E_{i-1}$ is semistable for $i \le k$. The quotient $V/ E_k$ need not be semistable, but \\
\\
$(\ast)$ for any isotropic subbundle $S \subset V$ containing $E_k$, we have $\mu (S/E_k) \le \mu(V/E_k)$.   \\
\\
It is easy to check that  $(E_i/E_{i-1}) \oplus (E_{i-1}^\perp / E_i^\perp)$ and $E_k^\perp / E_k$ have the induced  symplectic/orthgonal structure by  linear algebra.
 The semistability of $(E_i/E_{i-1}) \oplus (E_{i-1}^\perp / E_i^\perp)$ follows from 
 $E_{i-1}^\perp / E_i^\perp \ \cong \ (E_i/E_{i-1})^*$. 
 
Now we show that  $E_k^\perp / E_k$ is semistable. Since it is a symplectic/orthogonal bundle, it suffices to check that there is no destabilizing isotropic subbundle. Let $\widetilde{S}$ be  any isotropic subbundle of $E_k^\perp / E_k$. Then  $\widetilde{S} = S/E_k$, where $S$ is an isotropic subbundle of $E_k^\perp$ 
containing $E_k$. To prove the claim that $\mu (E_k^\perp /E_k) \ge \mu(\widetilde{S}) $, we note
\[
\mu (E_k^\perp /E_k) - \mu(\widetilde{S}) = ( \mu (E_k^\perp /E_k) - \mu (S) ) + (\mu (S) - \mu(\widetilde{S}) ).
\]
Since the above process ends up at $E_k$, we have $\mu(V) \ge \mu (S)$. 
From the isomorphism $V/E_k \cong (E_k)^\perp \otimes L$, we compute $\mu (E_k^\perp / E_k) = \mu(V)$. This shows that the first term on the righthand side is nonnegative. Also  the sign of the second term is same as that of $\mu (E_k) - \mu(S)$, which is nonnegative from the choice of $E_k$ with maximal slope. This shows the claim.
\end{proof}

Note that the uniqueness of the Harder--Narasimhan filtration of the underlying vector bundle guarantees the uniqueness of the sequence (\ref{seq}) for a symplectic/orthogonal bundle $V$.

In particular for $g=0$, the existence of a symplectic or orthogonal structure on a vector bundle $V$ on $\pp^1$ forces its splitting type to be  special. First note that the   form $\omega \colon V \otimes V \to \cO (\ell)$ can be normalized: 
when $\ell$ is even (resp. odd),  we may assume that $\ell=0$ (resp. $\ell=1$) by replacing $V$ by $V \otimes  \cO(\frac{\ell}{2})$ (resp.  $V \otimes  \cO(\frac{\ell-1}{2})$). 

\begin{cor} \label{HNP1}
Let $V$ be an $\cO(\ell)$-valued symplectic or orthogonal bundle over $\pp^1$ of the splitting type
\begin{equation*} \label{gensplit}
V \cong \cO(a_1)^{ r_1} \oplus \cO(a_2)^{ r_2} \oplus  \cdots  \oplus \cO(a_m)^{ r_m}
\end{equation*}
with $a_1>\cdots>a_m$. Then for  each $i $, we always have $ r_{m+1-i} =r_i $ and 
\begin{equation} 
 a_{m+1-i}  =-a_i \ \ \text{for} \ \ell=0 \ \ \text{and} \ \  a_{m+1-i}  =-a_i +1 \  \  \text{for}  \ \ell =1.
 \end{equation} 
 \end{cor}
 \begin{proof}
According to the type of Harder-Narasimhan filtration in  (\ref{HNvb}),   we have
\[
r_1 + (2n-r_m) = 2n, \ \ (r_1 + r_2) + (2n-(r_m+r_{m-1})) =2n, 
\]
and so on. This shows $ r_{m+1-i} =r_i $ for $1 \le i \le m$.

Also when $\ell=0$,  we have $\deg (V) = 0$ and $\deg (E^\perp) = \deg (E)$, so  the filtration (\ref{HNvb}) yields the equalities  
\[
a_1 r_1=  - a_m r_m,  \ \  a_1r_1 + a_2 r_2 =  - (a_m r_m + a_{m-1} r_{m-1}),
\]
and so on. This inductively shows $a_{m+1-i}  =-a_i $ for $1 \le i \le m$. 

When $\ell=1$, we have $\deg (V) = n$ and $\deg (E^\perp) = \deg (E) + n- \rk (E)$.   From this, the equality $a_{m+1-i}  =-a_i +1$  is obtained  similarly.
 \end{proof}

\section{Hecke curves} \label{Heckesection}
 \textbf{Throughout this section, we assume $n \ge 4$ for symplectic case and $n \ge 5$ for orthogonal case.} In \S 4.1, we establish a result on the uniqueness of a generic isomorphism as a technical toolkit.  And then we construct symplectic and orthogonal variants of Hecke curves  in \S 4.2 and 4.3 respectively.

\subsection{Uniqueness of a generic isomorphism} For  a nonnegative integer  $\delta$, a symplectic or orthogonal bundle $V$ is said to be \emph{$\delta$-stable} if 
\[
 \frac{\deg(E) +\delta}{\mathrm{rk} (E)} < \frac{\deg(V)}{n}   .   
 \]
 for each nonzero  isotropic subbundle $E$. This is a notion adapted to the symplectic and orthogonal case, which is originally the $(\delta, \delta)$-stability in \cite[Definition 5.1]{NR1}. Note that  the 0-stability is equivalent to the stability.  It is easy to see that the $\delta$-stability is an open condition in a family for each $\delta \ge 0$.

 \begin{lemma} \label{kkstable} 
 Let $\M$ denote either $\MS_C(n, L)$ or  an irreducible component of $\MO_C(n,L)$. Let $\delta$ be a positive integer.
 A general symplectic or orthogonal bundle in $\M$  is $\delta$-stable  if either  {\rm(}$g = \delta +2$ and $n > \delta +1${\rm)} or {\rm(}$g > \delta +2$ and $n > \frac{1}{2}(\delta +2)${\rm)}. 
\end{lemma}
 \begin{proof} A similar result for orthogonal case was proven in \cite[Proposition 3.5]{BG2}, but the above claim is slightly stronger.
 
 Let  $\mathcal{S}$ be the sublocus of $\M$ consisting of  bundles $V$   admitting an isotropic subbundle $E$ of rank $r$ and degree $d$ such that 
\begin{equation} \label{11ineqline}
\frac{d +\delta}{r} \ge \frac{\ell }{2 } \ , \ \ \text{or equivalently} \ \ d \ge \frac{1}{2} r \ell - \delta.
\end{equation}
It suffices to show that the dimension of  (any irreducible component of) $\mathcal{S}$  is strictly smaller than $ \dim (\M)$.
 
 From the exact sequence
\begin{equation} \label{extM}
0 \to E \to V \to (E^\perp)^*\otimes L \to 0,
\end{equation}
 we have
\begin{equation} \label{degEperp}
\deg (E^\perp) = d+(\frac{n}{2}-r)\ell.
\end{equation}
The inherited  form on the quotient $E^\perp /E $ is non-degenerate, so it is an $L$-valued symplectic/orthogonal bundle of rank $n-2r$. 
The above sequence can be put into the middle sequence of the following diagram:

\begin{equation} \label{symmdiagram}
\begin{tikzcd}[ arrows={-stealth}] 
 &  \null & E^*  L \rar[=] & E^*  L \\
0 \rar & E  \rar & V \rar \uar \ular[dash, dashed, red] & (E^\perp)^*  L \rar \uar & 0\\
 0 \rar & E  \rar \uar & E^\perp  \rar \uar & E^\perp/ E \rar \uar \ular[dash, dashed, red] & 0.
\end{tikzcd} 
\end{equation}

Note that this diagram is symmetric with respect to the diagonal  line (possibly with a change of sign). That is, the central term has the symmetry $V  \cong V^*  \otimes L$ and the right vertical sequence gives the bottom horizontal sequence after taking dual and tensoring $L$.  

The middle extensions are parametrized by 
$\pp H^1(C, \mathrm{Hom}(E^\perp \otimes E \otimes  L^*))$, which  appears again in the long exact sequence associated to the bottom sequence (tensored by $ E \otimes L^*$): 
\begin{equation} \label{extH1}
\to   H^1(E \otimes E \otimes  L^*) \to H^1(E^\perp \otimes E \otimes L^*) \to H^1((E^\perp/E) \otimes E \otimes  L^*) \to 0.
\end{equation}

Now we compute an upper bound on the dimension of those bundles $V$ fitting into  this diagram. 
A relevant result for this dimension count is in \cite[Proposition 2.6]{NR2}, which guarantees that for any given  family of vector bundles,
there is an  irreducible family  of  bundles whose generic member is stable, and  contains all the members of the given family. Since the  arguments in  \cite[Proposition 2.6]{NR2} can be adapted to the symplectic or orthogonal cases, we may assume  the stability of the objects   in the dimension count below.

Consider the diagram (\ref{symmdiagram}). For the right vertical sequence, 
\begin{itemize} 
\item $E$ moves in the moduli of vector bundles of rank $r$, which has  dimension $r^2(g-1)+1$.
\item $E^\perp / E$ moves in a family of  symplectic/orthogonal bundles of rank $n-2r$, which has dimension 
$ \frac{1}{2}(n-2r)(n-2r \pm 1) (g-1)$, where the sign $\pm$ refers to $+$ in the symplectic case and $-$ in the orthogonal case (throughout the proof). 
\item The bundle $F:= (E^\perp)^* \otimes L $ is obtained as an extension of $E^* \otimes L$ by $E^\perp / E $, which lies on $ \pp H^1 (C, (E^\perp / E )  \otimes E \otimes L^*)$. To count the dimension of deformations, we can assume that $E^\perp$ is stable, and so
\[
H^0(C, \mathrm{Hom}(E \otimes L^*,  E^\perp / E) )= 0.
\]
 Since $(E^\perp / E )  \otimes E \otimes L^*$ has rank $r(n-2r)$ and degree  $(n-2r)(d - \frac{1}{2}r \ell)$,  we have
 \[
 \dim \pp H^1 ( (E^\perp / E)  \otimes E \otimes L^*) \ = \    (n-2r)( \frac{1}{2}r\ell -d ) + r(n-2r)(g-1) -1.
 \]
\end{itemize}
These information determine all the terms in the diagram by symmetry mentioned above, except the central term $V$.   This is  determined by the choice of a point in $\pp H^1( C, E^\perp \otimes E \otimes L^*)$, which is in the middle of (\ref{extH1}).  Since a class in $H^1( (E^\perp/E) \otimes E \otimes L^*)$ is already chosen, it suffices to choose a class in $H^1( E \otimes E \otimes L^*)$.
By the (anti)-symmetry of the diagram \eqref{symmdiagram}, the extension class  lies on $H^1(\Sym^2 E \otimes L^*)$ and $ H^1( \wedge^2 E \otimes L^*)$ for the symplectic and orthogonal case respectively (see \cite[Criterion 2.1]{Hit}), which has dimension 
$(r \pm 1) (\frac{1}{2}r \ell  - d) + \frac{1}{2}r(r \pm 1) (g-1)$.

The total sum of the above moduli gives the upper bound:
\[
\dim (\mathcal{S}) \ \le \ (n-r \pm 1)(\frac{1}{2}r \ell -d) +   \frac{1}{2} \left( n^2 -2rn +  3r^2  \pm (n-r) \right) (g-1)  .
\]
Since $\dim (\M) = \frac{1}{2} n(n \pm 1) (g-1)$, the codimension of $\mathcal{S}$ in $\M$  is positive if
 \begin{equation} \label{codimSinM}
(n-r \pm 1)(\frac{1}{2}r \ell -d)  < \frac{r}{2}(2n-3r \pm 1)(g-1).
 \end{equation}
 When $g = \delta +2 $, this holds  by the inequality    \eqref{11ineqline} if 
\[
 (n-r \pm 1) \delta < \frac{r}{2}(2n-3r \pm 1)(\delta+1).
\]
This holds for $r=1$ if  $n > \delta +1$. Also if $2 \le r \le \frac{n}{2}$, then $n-r \pm 1 \le 2n-3r \pm 1$ from which  the above inequality follows. 
 Similarly, \eqref{codimSinM}  holds if $g > \delta +2$ and $n > \frac{1}{2}(\delta +2)$.
 This shows that the dimension of the locus $\mathcal{S}$ is strictly smaller than $\dim (\M)$ under the assumptions on $g$ and $n$.
 \end{proof}

\begin{lemma}  \label{generic1}
Let $\M$ denote either $\MS_C(n, L)$ or  an irreducible component of $\MO_C(n,L)$. 
Suppose that $[V]$ is a general point of  $\M$ representing a symplectic or orthogonal bundle $V$. Let $\phi \colon W \to V$ be a generic isomorphism of vector bundles with  $\deg (V) - \deg (W) = \delta >0$. Then $\dim H^0(C, W^* \otimes V) = 1$ if 
either {\rm(}$\delta =1$   and $ g \ge 3${\rm)}
or
\begin{equation} \label{genusbound}
g > \frac{3(\delta-1)n}{n-1}+1,  \ n > \frac{1}{2}(\delta+2).
\end{equation}
\end{lemma}
\begin{remark}  \label{ddstable}
(1) The argument below  is a refinement of the proof of \cite[Lemma 5.6]{NR1}, which was originally  given for vector bundles.  \\
(2) Note that in each case the condition on $g$ guarantees that $V$ is $\delta$-stable by Lemma \ref{kkstable}. It can be checked that for $\delta \ge 2$ and $n \ge 4$,
\[
\delta +2  < \frac{3 (\delta -1)n}{n-1} +1 .
\]
(3) Later we will apply this result for $\delta \le 4$, in which case the condition $n > \frac{1}{2}(\delta +2)$ can be removed.
\end{remark}
\begin{proof}[Proof of Lemma \ref{generic1}]
First we show that if $\dim H^0(C, W^* \otimes V)>1$, then there is a nonzero map $W \to V$ with a nontrivial kernel.

Choose a general $y \in C$. Then $ y \notin \mathrm{Supp}(D)$, where $D$ is an effective divisor  of degree $\delta <g$ such that 
$(\det V) \otimes (\det W^*) \cong \cO_C(D)$. Let $\psi \colon W \to V$ be a map which is not a constant multiple of $\phi$.  Since $\phi$ and $\psi$ are linearly independent at a general point $y$,  some linear combination $a\phi + b \psi$ is degenerate at $y$.  This map   cannot be a generic isomorphism: if it were,   the  quotient of $V$ by $(a \phi + b \psi) (W)$ is a torsion sheaf supported on a divisor containing $y$, which is a contradiction. This shows that $a \phi + b \psi \colon W \to V$ is a nonzero map with a nontrivial kernel.

Therefore, it suffices to show that those bundles $V$ which admit  a nonzero map  $\psi \colon W \to V$ with a nontrivial kernel are contained in a proper closed subset of $\M$.

Now let $V$ be a bundle in $\M$ admitting  a nonzero map  $\psi \colon W \to V$ with a nontrivial kernel $S$ and the quotient $E = W/S$ which is a subsheaf of $V$ of rank $r<n$. 
Let  $M:= E \cap E^\perp$. If $M$ is nonzero of rank $r_M>0$, then from the $\delta$-stability (see Remark \ref{ddstable} (2)), we have  
\[
\frac{\deg (M) +  \delta}{r_M} < \frac{\ell}{2}.
\]
Hence  as before in \eqref{degM}, 
\[
\deg (E) = \deg (M) + (r-r_M) \cdot \frac{\ell}{2}  < \frac{r \ell}{2} -  \delta.
\]
Also, since $S$ is a subsheaf of $V$ via $\phi$, we have $\deg (S) < \frac{1}{2}(n-r)\ell$. Therefore,
\[
\deg (W) = \deg S + \deg E < \frac{n \ell}{2} - \delta = \deg (W),
\]
which is a contradiction. 

This shows that $E \cap E^\perp =0$ and hence the restriction of the form $(V, \omega)$ to $E$ is non-degenerate  on a general   fiber.  Note that
\begin{equation} \label{degreeboundE}
\deg (E) = \deg (W) - \deg (S) > \frac{n \ell}{2} - \delta - \frac{(n-r)\ell}{2}   =   \frac{r \ell}{2} - \delta.
\end{equation}
Since 
\[
\deg (E^\perp) =  \deg (E) + \frac{(n-2r) \ell}{2} > \frac{(n-r)\ell}{2} - \delta,
\]
we have
\[
\deg (E) + \deg (E^\perp) \ge \deg (V) - 2 (\delta-1).
\]
 Hence we get a sequence
\begin{equation} \label{eeperptau}
0 \to E \oplus E^\perp \to V \to { \tau_0} \to 0,
\end{equation}
where ${ \tau_0} $ is a torsion sheaf of degree $ \le 2(\delta -1)$. 

Now we count dimension of bundles $V \in \M$ fitting into \eqref{eeperptau}.   
First, we note that the restriction of the form $\omega$ to $E$ yields a generic isomorphism $E \to E^* \otimes L$ such that the quotient $(E^* \otimes L) / E$ is a  torsion sheaf of degree $ -2 \deg (E) + r \ell <\delta$ by \eqref{degreeboundE}. Then $E$  is a subsheaf (or a  limit of subsheaves) of  an $L$-valued symplectic/orthogonal bundle $\widetilde{E}$ of rank $r$ with $\deg (\widetilde{E}/E)  = \frac{1}{2} r \ell - \deg (E) < \delta$.  Hence we get
\begin{eqnarray*}
\dim \{ E\} &\le& \dim \M + (\delta-1) \dim \pp (\widetilde{E}) \\
&=& \frac{r(r \pm 1)}{2} (g-1) +(\delta-1) r.
\end{eqnarray*}
(For the sign $\pm$, we take $+$ and $-$ sign in the symplectic and orthogonal case, respectively. The space $\M$ should refer to the corresponding space.)
Similarly, we have 
\[
\dim \{ E^\perp \} \ \le  \frac{(n-r)(n-r \pm 1)}{2} (g-1) +  (\delta-1) (n-r).
\]
Given two bundles $E$ and $E^\perp$, those bundles $V$ fitting into the sequence of the form \eqref{eeperptau} have dimension bounded from above by 
\[
2 (\delta -1) \cdot \dim \pp (E \oplus E^\perp) = 2 (\delta -1) n.
\]
Summing up these three dimensions, we get 
\[
\frac{1}{2}(n^2 - 2rn +2r^2 \pm n) (g-1) +3(\delta -1) n =: D_0. 
\]
From the inequality 
$D_0 \ < \ \dim (\M) = \frac{1}{2}n(n \pm 1)(g-1)$,
we get 
\[
(rn-r^2)(g-1) >3(\delta-1)n.
\]
The lefthand side is a quadratic polynomial in $r$, which has  minimum at $r=1$. Thus the inequality reduces to: $(n-1)(g-1) > 3 (\delta -1)n$. 
Hence it can be seen that $\dim (\M) > D_0 $ under the assumption that $g > \frac{ 3(\delta -1)n}{n-1} +1$. 

We conclude that a general bundle $V$ in $\M$ does not admit  a generic isomorphism $\phi \colon W \to V$  together with a map $\psi \colon W \to V$ having a nontrivial kernel.  Therefore, given a general $[V] \in \M$, a generic isomorphism $\phi \colon W \to V$ is  unique  up to a constant multiplication, if it exists.
\end{proof}

\subsection{Symplectic Hecke curves on $\M = \MS_C(n,L)$}\label{Heckecurve} \ 
 Let $(V, \omega)$ be a point of $\M$ representing an $L$-valued symplectic bundle $V$ of rank $n$. Now we construct  symplectic Hecke curves, which will turn out to have minimal degree on $\M$. The construction closely follows that of Hecke curves on $SU_C(n,d)$ in \cite{Hw} (and \cite{HR}). To construct a symplectic version of Hecke curves, we need to keep track of the symplectic forms in the process.

 Given  a linear functional $\theta  \in V^*|_x$ for some $x \in C$, let $V^\theta $ be the kernel sheaf of the composition $V  \to V|_x \stackrel{\theta }{\to} \cc$. Then we have an extension
\begin{equation} \label{Heckemu}
0 \to V^\theta  \to V \to (V|_x/\ker \theta) \otimes \cc_x \to 0,
\end{equation}
where $\cc_x$ is the skyscraper sheaf supported at $x$. 
  Taking dual of this Hecke transformation, we get
\[
0 \to V^* \to (V^\theta)^* \to \cc_x \to 0,
\]
which   restricts to the fiber at $x$ as an exact sequence of vector spaces:
\begin{equation} \label{Heckeatx}
0 \to   \langle \theta \rangle  \to V^* |_x \to (V^\theta)^*|_x \to  \cc \to 0.
\end{equation} 
Then there is a sheaf injection from $(V^\theta)^* (-x) := (V^\theta)^* \otimes \cO_C(-x)$ into $V^*$.
Consider the composition map:
\[
(V^\theta)^* (-x) \stackrel{\alpha}{\longrightarrow}  V^* \xrightarrow{\tom} V \otimes L^* \stackrel{\beta}{\longrightarrow} V^\theta(x) \otimes L^*,
\]
 where $\alpha$ and $\beta$ are dual to each other and $\tom$ is the skew-symmetric isomorphism associated to $\omega$. This gives a  skew-symmetric map
$ (V^\theta)^* \to V^\theta \otimes L^*(2x)$, or equivalently an $L^*(2x)$-valued form
\[
\textstyle{\bigwedge^2}(V^\theta)^* \to L^*(2x).
\] 
The restriction of this  form to the fiber $\bigwedge^2 (V^\theta)^*|_x  $ factors through $\bigwedge^2(\mathrm{Im} (\alpha_x))$, which is zero since  $\dim \mathrm{Im} (\alpha_x) = 1$.  Hence we get an $L^*(x)$-valued  skew-symmetric  form on $(V^\theta)^*$, or equivalently a skew-symmetric map
\[
\omega^\theta : (V^\theta)^* \to V^\theta \otimes L^* (x).
\]
 By computing the difference of degrees, we see that the subspace $\ker(\omega^\theta_x)  $  has codimension two in $ (V^\theta)^*|_x$.

For each linear functional  $\lambda$ of $(V^\theta)^*|_x$, let $\tV^\lambda$ be the bundle obtained by   the Hecke transformation
\[
0 \to \tV^\lambda \to (V^\theta)^* \to \left( (V^\theta)^*|_x / \ker(\lambda)\right) \otimes \cc_x \to 0.
\] 
\begin{lemma} \label{vanish}
The bundle  $\tV^\lambda $ is equipped with an $L^*$-valued   symplectic form  $\tom^\lambda $ induced from $\omega^\theta$ if and only if $ \ker (\omega^\theta_x) \subset \ker (\lambda)$. 
\end{lemma}
\begin{proof}
Consider the composition  map
\[
  \tV^\lambda  \stackrel{\ta}{\longrightarrow} (V^\theta)^*  \xrightarrow{\omega^\theta} V^\theta \otimes  L^*(x) \stackrel{\tb}{\longrightarrow} (\tV^\lambda)^* \otimes L^*(x),
\]
where $\ta$ and $\tb$ are dual to each other. This gives a form 
\[
\tom^\lambda: {\textstyle \bigwedge^2} \tV^\lambda \to L^*(x),
\]
which is nondegenerate outside the fiber at $x$. This factors through $\bigwedge^2 \mathrm{Im}(\ta)$, where $  \mathrm{Im}(\ta_x) = \ker (\lambda)$ has codimension one in $(V^\theta)^*|_x$. If $ \ker (\omega^\theta_x) \subset \ker (\lambda)$, then $\tom^\lambda_x$ further factors through the quotient space $\bigwedge^2  {\big(} (\ker(\lambda)/ \ker(\omega^\theta_x) {\big)}$, which is zero.  Therefore in this case, we have an $L^*$-valued form
\[
\tom^\lambda: {\textstyle \bigwedge^2} \tV^\lambda \to L^*.
\]
This is nondegenerate everywhere since $\tV^\lambda$ and $(\tV^\lambda)^* \otimes L^*$ have the same degree.

If $ \ker (\omega^\theta_x) \not\subset \ker (\lambda)$, it is easy to see that the  form $\tom^\lambda: {\textstyle \bigwedge^2} \tV^\lambda \to L^*(x)$ does not identically vanish on the fiber at $x$, and it remains to be a partially degenerate $L^*(x)$-valued form.
\end{proof}
Hence for each linear functional $ \lambda$ of $(V^\theta)^*|_x$ whose kernel contains $\ker(\omega^\theta_x)$,   we get a   bundle $(\tV^\lambda)^*$ equipped with an $L$-valued symplectic form $(\tom^\lambda)^*$.  Since this is parameterized by $\lambda \in \pp \big( (V^\theta)^*|_x / \ker(\omega^\theta_x) \big) \cong \pp^1$,  we get a family  
$\C^\theta[V]$
of $L$-valued symplectic bundles.
 In particular when $\ker(\lambda)$ coincides with the subspace $V^*|_x / \langle \theta \rangle$ from the sequence (\ref{Heckeatx}),   we get back $  (V, \omega) \cong \big( (\tV^\lambda)^*,   (\tom^\lambda)^* \big)$.

 This shows that  for each $\theta \in \pp V^*$, the family $\C^\theta[V]$   gives  a rational curve on $\M$ passing through the point  $[V]$, {provided that} it  is a nonconstant family and  the symplectic bundles $\tV^\lambda$'s are stable. 

\begin{lemma} \label{sympHecke}
 Assume $g \ge 3$ and  $n \ge 4$.
 Let $V $ be a general point of $\M$. Then 
 \begin{enumerate}
\item all the  symplectic bundles $\tV^\lambda$ appearing in $\C^\theta[V]$ are stable, and
\item  the map $\pp^1 \to \C^\theta[V]$ given by $\lambda  \mapsto \big( (\tV^\lambda)^*, (\tom^\lambda)^* \big) $ is generically injective.
 \end{enumerate}
 \end{lemma}
 \begin{proof}
 By Lemma \ref{kkstable}, we may assume that $V$ is 1-stable. 
 
 To show (1),  for each $[\tV^\lambda] $, let $E$ be an isotropic subbundle of $\tV^\lambda$. Then by the construction, $V^*$ has a subsheaf $\tE$ defined by the kernel of the composition $(E \to \tV^\lambda  \to (V^\theta)^*  \to  (V^\theta)^* / V^* \cong \cc_x)$. Then in any case $\deg (\tE) \ge \deg (E) -1$. By the $1$-stability of $V^*$, we have
 \[
 \frac{\deg(E)}{\rank (E)}  \le \frac{\deg(\tE)+1}{\rank (\tE)}  < \frac{ \deg(V^*)}{\rank (V^*)} = \frac{ \deg(\tV^\lambda)}{\rank (\tV^\lambda)}.
 \] 
 
 To show (2),  first we note that  a general member $\tV^\lambda$ is 1-stable, since the 1-stability is an open condition and the family contains a 1-stable bundle $V^*$. Suppose
 that  for $\lambda_1, \lambda_2 \in \pp^1$, both $\tV^{\lambda_1}$ and $\tV^{\lambda_2}$  are 1-stable and $\tV^{\lambda_1} \cong \tV^{\lambda_2}$. Then there are two linearly independent generic isomorphisms between $(V^\theta)^*$ and $\tV^{\lambda_1}$. By Lemma \ref{generic1}, this implies $\lambda_1 = \lambda_2$.  Therefore, a general member $\tV^\mu$ is not isomorphic to any other member.
   \end{proof}
 This rational curve $\C^\theta[V]$ is called the \textit{symplectic Hecke curve} on $\M$ associated to $\theta \in \pp V$. It will be shown in \S\:\ref{mainsection} that $\C^\theta[V]$ has degree $2n$.

\subsection{Orthogonal Hecke curves on $\MO_C(n,L)$} \label{Heckecurveortho} \  The construction  in the orthogonal case is similar to the symplectic case, but also there are some important differences. 
We follow the construction in \S\:\ref{Heckecurve}, pointing out the required modifications.

 Let $\M$ be any irreducible component of $\MO_C(n,L)$. Let $(V, \omega)$ be a point of $\M$ representing an $L$-valued  orthogonal  bundle $V$ of rank $n$. For $x \in C$, let $IG(2, V|_x)$ be the Grassmannian of 2-dimensional isotropic subspaces of $V|_x$.  Then for each $\Theta \in IG(2, V|_x)$, we can associate  a codimension 2 subspace $\Theta^\perp$ of $V|_x$.  Let $V^\Theta$ be the kernel sheaf of the composition $V \to V|_x \to (V|_x) / (\Theta^\perp)$. This gives an extension
\begin{equation} \label{Heckeext}
0 \to V^\Theta \to V \to (V|_x / \Theta^\perp)  \otimes \cc_x \to 0,
\end{equation}
Taking dual of this Hecke transformation, we get 
\[
0 \to V^* \to (V^\Theta)^* \to \cc_x^{\oplus 2} \to 0,
\]
 which restricts to the fiber at $x$ as an exact sequence of vector spaces
 \begin{equation} \label{HeckeTheta}
 0 \to \Theta \to V^*|_x \to (V^\Theta)^*|_x \to \cc^2 \to 0,
 \end{equation}
 since there is a canonical isomorphism $(V|_x / \Theta^\perp)^* \cong \Theta$.
 
Then  there is a sheaf injection from $ (V^\Theta)^*(-x) $ into $V^*$. Consider the composition map 
\[
(V^\Theta)^*(-x) \stackrel{\alpha}{\longrightarrow} V^*  \stackrel{\tom}{\to} V \otimes L^* \stackrel{\beta}{\longrightarrow} V^\Theta(x) \otimes L^*,
\]
where $\tom$ is the symmetric isomorphism associated to $\omega$. This gives a symmetric form 
\[
\Sym^2 (V^\Theta)^* \to L^*(2x).
\]
The restriction of this form to the fiber at $x$ factors through 
$\Sym^2 (\mathrm{Im} (\alpha_x)) = \Sym^2 \Theta$, which is zero since $\Theta$ is chosen to be isotropic. Hence we get an $L^*(x)$-valued symmetric form on  $(V^\Theta)^*$, or equivalently a symmetric map 
\[
\omega^\Theta : (V^\Theta)^* \to V^\Theta \otimes L^*(x).
\]
 By computing the difference of degrees, we see that the subspace $\ker(\omega^\Theta_x)$ has codimension four in $(V^\Theta)^*|_x$. Then a     subspace $\Lambda \subset (V^\Theta)^*|_x$ of codimension two which is isotropic with respect to this form corresponds to a 2-dimensional  isotropic subspace of the orthogonal vector space $(V^\Theta)^* / \ker(\omega^\Theta_x) \cong \cc^4$, and vice versa.

 For any subspace $\Lambda \cong \cc^{n-2}$ of $(V^\Theta)^*|_x$, let $\tV^\Lambda$ be the bundle obtained by the Hecke transformation
\[
0 \to \tV^\Lambda \to  (V^\Theta)^*  \to \big( (V^\Theta)^*|_x / \Lambda \big) \otimes \cc_x \to 0.
\]

\begin{lemma} 
The bundle  $ \tV^\Lambda $ is equipped with an $L^*$-valued   orthogonal form  $\tom^\Lambda $ induced from $\omega^\Theta$ if and only if and $\Lambda$ is an isotropic subspace containing $\ker (\omega^\Theta|_x)$. 
\end{lemma}
\begin{proof}
We can argue in the same way as in the proof of Lemma \ref{vanish}.  Consider the composition map 
\[
\tV^\Lambda \stackrel{\ta}{\longrightarrow} (V^\Theta)^*  \stackrel{\tom^\Theta}{\to} V^\Theta \otimes L^*(x) \stackrel{\tb}{\longrightarrow} (\tV^\Lambda)^* \otimes L^*(x),
\]
where $\ta$ and $\tb$ are dual to each other. This gives a form 
\[
\tom^\Lambda : \Sym^2 \tV^\Lambda \to L^*(x).
\]
This factors through $\Sym^2 (\mathrm{Im}(\ta))$, where $\mathrm{Im}(\alpha_x) = \Lambda$ has codimension two in $(V^\Theta)^*|_x)$. If $\Lambda \supset \ker (\omega^\Theta|_x)$, then $\tom^\Lambda_x$ further factors through the quotient space $\Sym^2 \big(   \Lambda / \ker (\omega^\Theta|_x) \big)$.  Now the form $\tom^\Lambda$ vanishes identically along the fiber at $x$  if $\Lambda$ is isotropic. Conversely,  $\tom^\Lambda$ does not vanish identically along the fiber at $x$  if either $\Lambda \not\supset \ker (\omega^\Theta|_x)$ or $\Lambda$ is not isotropic.
\end{proof}

The parameter space of codimension two isotropic subspaces containing  $\ker (\omega^\Theta|_x)$ is isomorphic to the isotropic Grassmannian $IG( 2, \big((V^\Theta)^*|_x /   \ker (\omega^\Theta|_x) \big) = IG(2,4)$. 
It is well-known that $IG(2, 4)$ is a disjoint union of two $\pp^1$'s.  For each $t \in IG(2, 4)$ representing such a codimension two  isotropic subspace $\Lambda_t$, let $\tV_t$ be the bundle $ (\tV^{\Lambda_t})^* $ equipped with the $L$-valued orthogonal form $\tom_t : = (\tom^{\Lambda_t})^*$.  Then we get a family 
\[
 \{ (\tV_t, \tom_t) \: : \: t \in IG(2, 4) \}
 \]
 of  $L$-valued orthogonal bundles.  In particular when $t=t_0$ corresponds to $\Lambda = V^*|_x / \Theta$  in (\ref{HeckeTheta}), 
   we get back the original bundle $ (\tV_{t_0}, \tom_{t_0} ) = (V, \omega)$. 
 
 Choose the $\pp^1$-component  of $IG(2, 4)$ containing $t_0$.\footnote{We expect that another component  (which does not contain $t_0$) produces orthogonal bundles with 2nd Stiefel--Whitney class different from $w_2(V)$.} Then  the family
 \[
\C^\Theta[V] \ := \ \{ (\tV_t, \tom_t)  \: : \: t \in \pp^1 \} 
 \]
 is a rational curve on $\M$ passing through the point  $[V]$, {provided that} the bundles $\tV_t$'s are stable as  orthogonal bundles. 

\begin{lemma} \label{orthoHecke}
 Assume $g \ge 5$ and $n \ge 5$.
 Let $V $ be a general point of $\M$. Then 
 \begin{enumerate}
\item all the  orthogonal bundles $\tV_t$ appearing in $\C^\Theta [V]$ are stable, and
\item  the map $\pp^1 \to \C^\Theta[V]$ given by $t \mapsto (\tV_t, \widetilde{\omega}_t) $ is generically injective.
 \end{enumerate}
 
 \end{lemma}
  
   \begin{proof}
 By Lemma \ref{kkstable}, we may assume that $V$ is 2-stable. 
 
 The stability of $\tV_t$ is equivalent to the that of its dual: $\tV^{\Lambda_t}$. 
Let $E$ be an isotropic subbundle of $\tV^{\Lambda_t}$. Then by the construction, $V^*$ has a subsheaf $\tE$ defined by the kernel of the composition $(E \to  \tV^{\Lambda_t} \to (V^\Theta)^* \to \cc_x^{\oplus 2})$. Then in any case $\deg (\tE) \ge \deg (E) -2$. By the stability condition of $V^*$, we have
 \[
 \frac{\deg(E)}{\rank (E)}  \le \frac{\deg(\tE)+2}{\rank (\tE)}  < \frac{ \deg(V^*)}{\rank (V^*)} =\frac{ \deg(\tV^{\Lambda_t})}{ \rank (\tV^{\Lambda_t})}.
 \] 
 
To show (2),   first we note that  a general member $\tV_t$ is 2-stable, since the 2-stability is an open condition and the family contains a 2-stable bundle $V$. Suppose  that  for $s, t \in \pp^1$, both $V_s$ and $V_t$  are 2-stable and $\tV_s \cong \tV_t$. Then there are two linearly independent generic isomorphisms between $V^\Theta$ and $\tV_t$. By Lemma \ref{generic1}, this implies $s=t$.  Therefore, a general member $\tV_t$ is not isomorphic to any other member, and we get the generic injectivity.
 \end{proof}

This rational curve $\C^\Theta[V]$ is called the \textit{orthogonal Hecke curve} on $\M$ associated to $\Theta \in IG(2, V|_x)$. It will be shown in \S\:\ref{mainsection} that $\C^\Theta[V]$ has degree $4n$.

\begin{remark}
In Lemma \ref{sympHecke} and \ref{orthoHecke}, we have seen that  the symplectic/orthogonal Hecke curve is  a generically injective image of  $\pp^1$. 
If we assume a higher genus bound, then we can guarantee the 2- or 4-stability of a general member of $\M$, and show the injectivity of the maps $\pp^1 \to \C^\theta[V]$ or $\pp^1 \to \C^\Theta[V]$.
\end{remark}

\section{minimal rational curves} \label{mainsection}
 
Based on  the results  in \S\:\ref{Harder} and \S\:\ref{Heckesection}, we can now compute the minimal degree of rational curves on the moduli spaces by adapting the computation in \cite[\S\:2]{Sun} to our situation. In doing that, we need to additionally keep track of the form $\omega$.

\subsection{Degree formula of rational curves} \label{degformula}
Let $\M$ be either  $\MS_C(n, L)$ or an irreducible component of $\MO_C(n,L)$. 
Let $u : \pp^1 \to \M$ be a  generically injective morphism.  
 Let $\M^{reg} \subset \M$  be the open sublocus consisting of those symplectic/orthogonal  bundles whose underlying vector bundle is stable. Let $\pi_1 \colon  C \times \pp^1 \to C$ and $\pi_2 \colon C \times \pp^1 \to \pp^1$ be the projections. Let $\cL :=   \pi_1^* L$.
\begin{lemma} \label{globalomega}
If $u (\pp^1) \subset \M^{reg}$,
 there is a vector bundle $\V$ over $C \times \pp^1$ equipped with a form $\Omega \colon  \V \otimes \V \to \cL  \otimes  \pi_2^*  \cO_{\pp^1}(\td)$ for some $\td$ such that for each $t \in \pp^1$, the restriction of $(\V, \: \Omega)$ to $\pi_2^{-1} (t) = C \times \{t\}$ is isomorphic to the symplectic/orthogonal bundle represented by $u( t) $. 
 \end{lemma}
 
 \begin{proof}
 Let $V_t$ be the  bundle represented by $u(t)$ for  $t \in \pp^1$. First  we note that $\dim H^0(C, \mathcal{E}nd (V_t, V_t^* \otimes L)) = 1$ for each $t$ since 
 $V_t \cong V_t^* \otimes L$  and 
$V_t$ is  a stable vector bundle. 
 
 By composing with the generically injective morphism $f: \M \to SU_C(n, \frac{1}{2}n\ell)$, the curve $(f \circ u)(\pp^1)$ lies in the moduli of stable vector bundles over $C$.
 By \cite[Lemma 2.1]{Sun},  there is a vector bundle  $\V$ over $C \times \pp^1$  such that for each $t \in \pp^1$, the restriction of $\V$ to $\pi_2^{-1} (t) = C \times \{t\}$ is isomorphic to the vector bundle  $V_t$.   From the above observation, we see that $(\pi_1)_*  \mathcal{E}nd (\V, \V^* \otimes \cL)$ is a line bundle over $\pp^1$, say $\cO_{\pp^1} (-\td)$ for some $\td$. By adjunction formula we have
 \[
 \Hom(\cO_{\pp^1} (-\td),(\pi_1)_*  \mathcal{E}nd (\V, \V^* \otimes \cL))\cong \Hom((\pi_1)^* \cO_{\pp^1} (-\td), \mathcal{E}nd (\V, \V^* \otimes \cL)),
 \]
so   there is a nonzero map 
 \[
\cO_{C\times\PP^1}   \rightarrow \mathcal{E}nd (\V, \V^* \otimes \cL) \otimes  \pi_1^*  \cO_{\pp^1} (\td) .
 \]
 This gives a map $\Omega: \V \to \V^* \otimes \cL \otimes \pi_1^*  \cO_{\pp^1} (\td) $ with the desired property.
 \end{proof}
Note that the   form $\Omega \colon \V \otimes \V \to  \pi_1^*  \cO_{\pp^1} (\td)$ can be normalized: 
When $\td$  is even (resp. odd),  we may assume that $\td=0$ (resp. $\td=1$) by replacing $\V$ by $\V \otimes   \pi_1^*  \cO_{\pp^1}(\frac{\td}{2})$ (resp.  $V \otimes  \pi_1^*   \cO_{\pp^1}(\frac{\td-1}{2})$). 
The normalized family $\V$ can be viewed as a family of $\cO(\td)$-valued symplectic/orthogonal bundles over $\pp^1$ parametrized by $C$ for $\td=0,1$. For  $x \in C$,  we write $\PP^1_{x} := \pi_1^{-1}(x)$.

Under this normalization,  the generic splitting type of $\V$ is given as in Corollary \ref{HNP1}.
 Then there is a relative Harder--Narasimhan filtration
\[
0 =\cE_0 \subset \cE_1  \subset \cdots  \subset \cE_{m} = \V
\]
which restricts to a generic fiber ${\pp^1_x}$ as $(\cE_i/\cE_{i-1})|_{\pp^1_x}  \cong \cO_{\pp^1}(a_i)^{\oplus r_i}$ for $1 \le i \le m$. If we let $\cF_i = \cE_i/\cE_{i-1}$, the bundle $\cF'_i := \cF_i \otimes \pi_2^* \cO_{\pp^1} (-a_i)$ is torsion-free with generic splitting type $\cO_{\pp^1}^{r_i}$. 

For a general $p \in \pp^1$,  the restriction of $\V$ to $\pi_2^{-1} (p) = C \times \{ p \}$ is denoted simply by $V$ and  we write $E_i := \cE_i|_{ C \times \{p \}}$, $F_i := \cF_i|_{ C \times \{p \}}$. 
By the degree formula  of \cite[(2.2)]{Sun}, we have
\begin{equation}\label{delta}
 \deg(u(\pp^1))= 2n \left( \sum\limits^{m}_{i=1}c_2(\cF_i') + \sum\limits_{i=1}^{m-1} (a_i-a_{i+1})(\mu(V)-\mu(E_i)) r_i \right).
\end{equation}
Furthermore, we observe:
\begin{itemize}
\item The bundles $E_i$ are isotropic and $E_i^\perp = E_{m-i}$  for $1 \le i \le \lceil \frac{m}{2} \rceil$.
\item Hence the summation for $\lceil \frac{m}{2} \rceil \le  i \le  m$  can be computed from the summation  for $1 \le i \le \lceil \frac{m}{2} \rceil$ by using the formula 
\[
(\mu(V)-\mu(E_i^{\perp})) \rk (E_i^{\perp}) \: = \: (\mu (V) - \mu(E_i) ) \rk (E_i),
\]
which comes from $\deg  (E_i^{\perp}) = \deg (E_i) + ( \frac{n}{2} - r_i ) \ell$.
\item In particular,  \eqref{delta} can be computed using   the isotropic bundles only.
\end{itemize}
 \subsection{Minimality of symplectic and orthogonal Hecke curves} First consider the moduli space $\M = \MS_C(n, L)$ of symplectic bundles.
\begin{theorem} \label{mainHecke}
 Assume $g \ge 3$ and  $n \ge 4$. Suppose that $u(\pp^1) \subset \M$ is a rational curve passing through a general point. Then $u(\pp^1)$ has degree $\ge  2n$.  Also if $u (\pp^1)$ has degree $2n$,  it is a symplectic Hecke curve.
\end{theorem}
\begin{proof}
We use the notations from  \S\:\ref{degformula}.
By \cite[Lemma 2.2]{Sun}, $c_2(\cF_i') \geq 0$ in   \eqref{delta} since the generic splitting type of $\cF_i'$ is trivial. By Lemma \ref{kkstable},   the curve $u \colon \pp^1 \to \cM$ passes through a point $[V]$ corresponding to a $1$-stable symplectic bundle. Hence we have $(\mu (V) - \mu(E_i) )r_i > 1$. These together show that $\deg(u(\pp^1)) \geq 2n$. 
 
 If the equality holds, then   $m=1$  and $c_2(\V) = 1$. In particular, $\td$ must be even. 
 Now we can apply the same argument as in \cite{Sun}:
The bundle $\V$ has exactly one jumping line
\[
\V_{\PP^1_{x}} \cong \cO(-1)\oplus \cO(1) \oplus \cO^{ n-2}
\]
(\cite[Lemma 2.5]{Sun}) and $\V$ admits  an elementary transformation
\[
0 \to \pi_1^*W \to \V \to \cO_{\PP^1_{x}}(-1) \to 0
\]
for some vector bundle $W$ on $C$ by \cite[Lemma 2.2]{Sun}.
Therefore $\cV$ induces a symplectic Hecke curve as was defined in \S\:\ref{Heckecurve}.  (In  \S\:\ref{Heckecurve}, the dual family $\V^*$ was constructed as a Hecke transformation of a fixed bundle and then its dual was taken at the final stage.)
\end{proof}

Now consider the moduli space  $\MO_C(n, L)$ of orthogonal bundles.  Let $\M$ be any irreducible component of $\MO_C(n, L)$. We note that the degree bound is twice as large as that of symplectic case.

\begin{theorem} \label{mainHeckeortho}
Assume   $g \ge 5$ for $n\ge5$. Suppose that $u(\pp^1) \subset \M $ is a rational curve passing through a general point. Then $u(\pp^1)$ has degree $\ge  4n$.  Also if  it has degree $4n$,  it is an orthogonal Hecke curve.
\end{theorem}
\begin{proof}
As before we have $c_2(\cF_i') \ge  0$. Since the curve $u \colon \pp^1 \to \cM$  passes through a point $[V]$ corresponding to a $2$-stable orthogonal bundle by  Lemma \ref{kkstable}, we have $(\mu (V) - \mu(E_i) )r_i > 2$. Therefore, $\deg(u(\pp^1)) \ge 2n$.
 (Also we observe that if $\deg(u(\pp^1)) \le 4n$, then  $m=1$  and $c_2(\V) = 1$. In particular, $\td$ must be even. )
 
 Now we exclude the possibility of degree $2n$.
  If  $\deg(u(\pp^1)) = 2n$, then
 again by applying   \cite[Lemma 2.2 and 2.5]{Sun},  we see that $\V$ admits an elementary transformation of the type:
 \begin{equation} \label{4ncase}
0 \to \pi_1^*W  \to \V \to \cO_{\PP^1_{x}}(-1) \to 0
\end{equation}
for some vector bundle $W$.  If this defines a 1-parameter family of orthogonal bundles, then each $\V|_{C \times \{ p \} }$ induces  a symmetric  form $\omega_p \colon \Sym^2(W) \to L(x)$ with a  kernel of codimension 2 in $W|_x$. Since $\ker (\omega_p)$ is fixed from the sequence \eqref{4ncase}, the form $\omega = \omega_p$ is independent of $p$. 
But  when we take elementary transformations associated to a one-dimensional subspaces in $W|_x$,  even if it is chosen as isotropic subspace, 
the  induced symmetric forms of the resulting bundles 
$\V|_{C \times \{ p \} }$   do not  vanish identically on the fiber at $x$ except two choices. (This is a consequence of the fact that the isotropic Grassmannian $IG(1,2)$ consists of two points.) Hence the family $\{ \V|_{ C \times \{ p \} } \: : \: p \in \pp^1 \}$ in \eqref{4ncase} do not lie inside $\MO_C(n,L)$.

So far we have shown that $\deg u( \pp^1) \ge 4n$. 
 Now  we characterize the curves of degree $4n$. If   $\deg(u(\pp^1)) = 4n$, there are three  possible types of the associated elementary transformations: Namely,
  \begin{equation} \label{8ncase-1-1}
0 \to \pi_1^*W  \to \V \to \cO_{\PP^1_{x_1}}(-1) \oplus \cO_{\PP^1_{x_2}}(-1)  \to 0
\end{equation}
for $x_1 \neq x_2 $ in $C$ or
\begin{equation} \label{8ncase-2}
0 \to \pi_1^*W \to \V \to \cO_{\PP^1_{x}}(-2) \to 0,
\end{equation}
or 
\begin{equation} \label{8ncase-1}
0 \to \pi_1^*W  \to \V \to \cO_{\PP^1_{x}}(-1)^{\oplus 2} \to 0.
\end{equation}
We now rule out the first two cases. The case  \eqref{8ncase-1-1} is ruled out  by the same reason  as \eqref{4ncase}:  since the elementary transformations  are taken at two fibers separately, we get only $2 \cdot 2 =4$ bundles in the family $\{ \V|_{ C \times \{ p \} } \: : \: p \in \pp^1 \}$ on which the induced  $L(x_1+x_2)$-valued  symmetric form vanishes identically on the fibers at $x_1,x_2$. 

The case \eqref{8ncase-2} can be understood as a consecutive  elementary transformations at the fiber $W|_x$ such that 
the first one is taken for a fixed choice of a point in $\pp(W|_x)$ and the second one is taken for a line in $\pp(W|_x)$. The same problem arise in this process as in \eqref{4ncase}. 

In the remaining case \eqref{8ncase-1}, the sequence restricts to each $C \times \{ p \}$  as an elementary transformation of the form
\[
0 \to W \to \V|_{C \times \{p\}} \to \cc_x^{\oplus 2} \to 0.
\]
 Then we can see that the rational curve $\{ \V|_{ C \times \{ p \} } \: : \: p \in \pp^1 \}$ lies on a component $\M$ of $\MO_C(n,L)$ only if it is an orthogonal Hecke curve as was defined in \S\:\ref{Heckecurveortho}. (In \S\:\ref{Heckecurveortho}, we constructed the dual family $\V^*$ first as a Hecke transformation of a fixed bundle and then took its dual at the final stage. This corresponds to the above construction in family.)
\end{proof}

\section{Applications} \label{applications}

To discuss applications of symplectic and orthogonal Hecke curves, we need to check that the symplectic (resp. orthogonal) Hecke curves $\C^\theta[V]$ (resp. $\C^\Lambda[\Theta]$) are effectively parameterized by $\theta \in \pp(V^*)$ (resp. by $\Theta \in IG(2, V)$). 

\begin{lemma}  \label{effpara}
\begin{enumerate}
\item For $n \ge 4$, assume  $g > 4+ \frac{3}{n-1}$ and let   $[V] \in \MS_C(n,L)$ be a  general point.  If $\theta_1 \neq \theta_2$ in $\pp(V^*)$, then  $\C^{\theta_1}[V]$ and $\C^{\theta_2}[V]$ are different symplectic Hecke curves.
\item For $n \ge 5$, assume    $g > 10 + \frac{9}{n-1}$ and let $[V] \in \MO_C(n,L)$ be a general point. If $\Theta_1 \neq \Theta_2$ in $IG(2,V)$, then $\C^{\Theta_1}[V]$ and $\C^{\Theta_2}[V]$ are different  orthogonal Hecke curves. 
\end{enumerate}
\end{lemma}

\begin{proof}
We exploit the notations from \S 4.
Assume $\theta_1 \neq \theta_2$. We show that a general  $[V_2] \in \C^{\theta_2}[V]$ is not isomorphic to  any $[V_1] \in \C^{\theta_1}[V]$.  Suppose  there were an isomorphism $\psi \colon V_1 \cong V_2$. Let $V^{\theta_1, \theta_2}$ be the subsheaf of $V$ given by $V^{\theta_1, \theta_2} = V^{\theta_1} \cap V^{\theta_2}$, so that $\deg ( V^{\theta_1, \theta_2} ) = \deg (V) -2$. Then there are two generic isomorphisms $V^{\theta_1, \theta_2} \to V_2$, one given by the composition
\[
V^{\theta_1, \theta_2} \subset V^{\theta_1} \subset V_1 \stackrel{\psi}{\to} V_2  
\]
and another by the composition
\[
V^{\theta_1, \theta_2} \subset V^{\theta_2} \subset V_2 . 
\]
Since $[V_2]$ is general in $\C^{\theta_2}[V]$ and $[V] \in \C^{\theta_2}[V]$, the point $[V_2]$  is as much general as $[V]$, and so $V_2$ is $2$-stable.  Thus we may apply Lemma \ref{generic1} to the map $V^{\theta_1, \theta_2} \to V_2$ in place of $W \to V$.
By the  assumption on $g$ and Lemma \ref{generic1} for $\delta=2$, these two generic isomorphisms coincide.  

Now consider the dual map $V_2^* \to (V^{\theta_1, \theta_2})^*$ of this generic isomorphism.
From the above observations, the restriction of this dual map  to the fiber at $x$  has the image inside the kernel of two maps for $i=1,2$:
\[
(V^{\theta_i})^*|_x \to  (V^{\theta_1, \theta_2})^*|_x .
\]
This shows $ V_2 \cong V$, which is a contradiction. 

The same argument shows  (2). In this case, we  apply Lemma \ref{generic1} for $\delta =3,4$ to show the coincidence  of two generic isomorphisms. (In modifying the argument, the genus assumption should be adapted correspondingly). 
\end{proof}

\begin{rmk}
To show that the spaces $\pp(V^*)$  and $IG(2, V)$ effectively parameterize the symplectic  and orthogonal Hecke curves respectively, one may instead try to show the injectivity of the tangent maps 
\[
\pp (V^*) \to \pp T_{[V]} \MS_C(n,L) \ \ \text{and} \ \ IG(2,V) \to \pp T_{[V]} \MO_C(n,L)
\]
sending a symplectic or orthogonal Hecke curve to its tangent at $[V]$. We believe this can be shown to be an embedding under certain assumptions on $g$ and $n$ as in \cite[theorem 3.1]{HR}, but we do not pursue it here. 
\end{rmk}

By Theorems \ref{mainHecke}, \ref{mainHeckeortho}, and Lemma \ref{effpara}, we get the following result.

\begin{cor} \label{Chow}
In the same notation as Lemma \ref{effpara}, the space $\pp (V^*)$ {\rm(}resp. $IG(2,V)${\rm)} is the normalization of the Chow variety of rational curves passing through $[V]$ in $\MS_C(n, L)$ {\rm(}resp. an irreducible component of $\MO_C(n,L)${\rm)}.
\end{cor}

Now we discuss the nonabelian Torelli theorem for the moduli of symplectic and orthogonal bundles.  This was proven before in  \cite[Theorem 0.1]{BH} for principal bundles in general, and later reproved in \cite[Theorem 4.3]{BGM} for the symplectic case. The approach in  \cite{BH} was to examine the strictly semistable locus and reduce to the classical Torelli theorem for the principally polarized Jacobians.  On the other hand,  both the proof of  \cite{BGM} and ours   use the Hecke correspondence, focusing on  the stable locus. But  we get the wanted result directly from the minimal rational curves on $\M$.

\begin{theorem} \label{Torelli}
Let $C_1$ and $C_2$ be smooth algebraic curves of genus $g$. Also let  $L_1$ and $L_2$ be line bundles of the same degree $\ell$ over $C_1$ and $C_2$, respectively. 
\begin{enumerate}
\item Assume $g > 4 + \frac{3}{n-1}$ for $n \ge 4$. If the moduli spaces $\MS_{C_1}(n,L_1)$ and $\MS_{C_2}(n,L_2)$ are isomorphic, then $C_1$ and $C_2$ are isomorphic.
 \item   Assume   $g > 10 + \frac{9}{n-1}$ for $n \ge 5$. Let $\M_1$ and $\M_2$ be any irreducible components of $\MO_{C_1}(n,L_1)$ and $\MO_{C_2} (n, L_2)$, respectively.  It $\M_1$ and $\M_2$ are isomorphic, then $C_1$ and $C_2$ are isomorphic. 
 \end{enumerate}
\end{theorem} 
\begin{proof} 
(1) \ Given an isomorphism $\Phi \colon \MS_{C_1}(n,L_1) \cong \MS_{C_2}(n,L_2)$, 
let $V_1$ be a symplectic bundle over $C$ corresponding to a general point $[V_1]$ of $\MS_{C_1}(n,L_1) $, and $V_2$ be the symplectic bundle over $C_2$ corresponding to the general point $\Phi([V_1])$ of $\MS_{C_2}(n,L_2)$. Then
$\Phi$ sends  the rational curves of degree $2n$ passing through $[V_1] \in \MS_{C_1}(n,L_1)$ to those  rational curves through $[V_2] \in \MS_{C_2}(n,L_2)$.  By Corollary \ref{Chow}, we get the induced isomorphism  $\pp (V_1^*) \cong \pp (V_2^*)$. From this isomorphism of  rational fibrations, we get the induced isomorphism  $C_1 \cong C_2$  on the base curves.\\
(2) \  In the same way, the isomorphism $\Phi \colon \M_1 \cong \M_2$ induces the isomorphism of isotropic Grassmannians $IG(2, V_1)$ and $IG(2, V_2)$, where $\Phi([V_1]) = [V_2]$. From this, we get an isomorphism $C_1 \cong C_2$.
\end{proof}

We can also classify the  automorphisms on the moduli spaces of symplectic or orthogonal bundles. The symplectic case has been shown  in \cite[Theorem 6.1]{BGM} by another method. Below we discuss the orthogonal case only. 
\begin{theorem}  \label{automorphism} 
Let $\M_0$ be the  component of $\MO_C(n, \cO_C)$ containing the trivial orthogonal bundle $\cO_C^{\oplus n}$. If $n \ge 5$ and $g > 10 + \frac{9}{n-1}$,  every automorphism of $\M_0$  is   induced by an automorphism of $C$ and a line bundle of order two.
\end{theorem}
\begin{proof}
Let $\sigma$ be an automorphism of $\M_0$. Since the degree of rational curves is fixed under $\sigma$,  we get an induced   isomorphism $\widetilde{\sigma}$ between isotropic Grassmannian bundles $\pi \colon IG(2,V) \to C$ and $\pi' \colon IG(2,V') \to C$, where $V'$ is the bundle represented by $\sigma([V])$.  Note that $\widetilde{\sigma}$ sends a fiber to a fiber, since $g \ge 2$. Thus by composing with an automorphism $\tau$ of $C$, we may assume that $\pi = \pi' \circ \widetilde{\sigma}$.

Let $\mathcal{U^*} \to IG(2,V)$ and $(\mathcal{U}')^* \to IG(2,V')$ be the dual of the relative universal bundles. Then  $\widetilde{\sigma}^* (\mathcal{U'})^* \cong \mathcal{U}^* \otimes \pi^* N$ for some line bundle $N$ over $C$. By taking push-forward, we have
\[
 \pi_* \widetilde{\sigma}^* (\mathcal{U'})^* = \pi'_* \widetilde{\sigma}_* \widetilde{\sigma}^* (\mathcal{U'})^* \cong V'
\]
and also
\[
\pi_* \widetilde{\sigma}^* (\mathcal{U'})^*  \cong   \pi_* (\mathcal{U}^* \otimes \pi^* N)  \cong V \otimes N.
\]
Hence $V' \cong  V \otimes N$ for  a line bundle $N$, up to  $\tau \in \mathrm{Aut}(C)$.   

Now it remains to show that $N$ is a line bundle of order two.
Since both $V$ and $V'$ are $\cO_C$-valued symplectic or orthogonal bundles, we get
\[
V' \cong V \otimes N \cong V^*  \otimes N \cong ((V')^* \otimes N)  \otimes N \cong V' \otimes N^2.
\]
This shows that the morphism $[V] \mapsto [V \otimes N^2]$ is the identity map  on $\M_0$. Hence we  can see that  $N^2 \cong \cO_C$ by specializing $V$ to the trivial orthogonal  bundle $V \cong \cO_C^{\oplus n} $.  
\end{proof}


\section*{acknowledgenent} The first named author would like to thank George H. Hitching for a helpful suggestion on the proof of Lemma \ref{generic1}. The first and the third named author thank Jaehyun Hong and Kyeong-Dong Park for  valuable discussions on isotropic Grassmannians. 

I. Choe was supported by Basic Science Research Programs through the National Research Foundation of Korea (NRF) funded by the Ministry of Education (NRF-2020R1F1A1A01068699).


\begin{thebibliography}{99}
\bibitem{BLS} Beauville, A.; Laszlo, Y.; Sorger, C.: \textsl{The Picard group of the moduli of $G$-bundles on a curve}, Compos. Math. \textbf{112} (1998), no. 2, 183--216.

\bibitem{BG1} Biswas, I.; G\'{o}mez, T. \ L.: \textsl{Hecke correspondence
for symplectic bundles with application to the Picard bundles},
Internat. J. of Math., \textbf{17} (2006), no. 1, 45--63.


\bibitem{BG2} Biswas, I.; G\'{o}mez, T. \ L.: \textsl{Hecke transformation for orthogonal bundles and stability of Picard bundles}, Comm. Anal. Geom. \textbf{18}  (2010), no. 5, 857--890.

\bibitem{BGM} Biswas, I.; G\'{o}mez, T. \ L.;  Mu\~{n}oz, V.: \textsl{Automorphisms of moduli spaces of symplectic bundles}, Internat. J. Math. \textbf{23} (2012), no. 5, 1250052, 27.

\bibitem{BH} Biswas, I.; Hoffmann, N.: \textsl{Torelli theorem for moduli spaces of principal bundles over a curve}, Ann. Inst. Fourier (Grenoble) \textbf{62} (2012), no. 1, 87--106.



\bibitem{Hit} Hitching, G. \ H.: \textsl{Subbundles of symplectic and orthogonal vector bundles over curves}, Math. Nachr. \textbf{280} (2007), no. 13-14, 1510--1517 

\bibitem{Hw1} Hwang, J.-M.: \textsl{Tangent vectors to Hecke curves on the moduli space of rank 2  bundles over an algebraic curve}, Duke J. Math. \textbf{101} (2000), 179--187.

\bibitem{Hw} Hwang, J.-M.: \textsl{Hecke curves on the moduli space of vector bundles over an algebraic curve}, Algebraic Geometry in East Asia, (Kyoto, 2001), 155--164.

\bibitem{HR} Hwang, J.-M.; Ramanan, S.: \textsl{Hecke curves and Hitching discriminant}, Ann. Sci. \'{E}cole Norm. Sup. (4) \textbf{37} (2004), no. 5, 801--817.



\bibitem{Mum} Mumford, D.:  \textsl{Theta characteristics of an algebraic curve}, 
Ann. Sci. \'{E}cole Norm. Sup. (4) \textbf{4} (1971), p. 181--192.

\bibitem{NR1} Narasimhan, M. \ S.; Ramanan, S.: {Geometry of Hecke cycles I},  In: C. P. Ramanujam --a tribute, Springer Verlag (1978), 291--345.

\bibitem{NR2} Narasimhan, M. \ S.; Ramanan, S.: \textsl{Deformation of the moduli space of vector bundles over an algebraic curve}, Ann. of Math. (2) \textbf{101} (1975),  391--417.

\bibitem{Ra} Ramanan, S.: \textsl{Orthogonal and spin bundles over hyperelliptic curves}, Proc. Indian Acad. Sci. Math. Sci. \textbf{90} (1981), no. 2, 151--166.

\bibitem{Ramanathan} Ramanathan, A.: \textsl{Moduli of principal bundles over algebraic curves}: I and II, Proc. Indian Acad. Sci. (Math. Sci.) \textbf{106} (1996), 301--328, 421--449.

\bibitem{Ser} Serman, O.: \textsl{Moduli spaces of orthogonal and symplectic bundles over an algebraic curve}, Compos. Math. \textbf{144} (2008), no. 3, 721--733.

\bibitem{Ser2} Serman, O.: \textsl{Moduli spaces of orthogonal bundles over an algebraic curve}, Th\'{e}se, Universit\'{e} Nice Sophia Antipolis, (2007), https://tel.archives-ouvertes.fr/tel-00262100.

\bibitem{Sun} Sun, X.: \textsl{Minimal rational curves on moduli spaces of stable bundles}, Math.  Ann.  \textbf{331} (2005), no. 4, 925--937.

\end{thebibliography}
\end{document}